\documentclass[12pt]{article}
\usepackage{hyperref}
\usepackage{amssymb}
\usepackage{amsthm}
\usepackage{amsmath} 
\usepackage{xcolor}

\usepackage{graphicx}

\usepackage{tikz}
\usepackage{pgfplots}
\usetikzlibrary{patterns}
\usetikzlibrary{perspective}
\usetikzlibrary{arrows}
\usetikzlibrary{3d,calc}

\usepackage{stmaryrd,bm}
\usepackage{accents}
\usepackage[normalem]{ulem} 
\usepackage{mathrsfs}
\usepackage{enumitem}
\usepackage{isomath}
\usepackage{rotating}
\usepackage{cite}
\usepackage{datetime}

\usepackage{tikzscale}
\pgfplotsset{compat=1.18}
\usetikzlibrary{angles, calc, intersections, quotes, arrows.meta}
\usepackage{tkz-euclide}
\usepackage{color}
\usepackage{soul}
\newcommand{\dst}{\displaystyle}
\renewcommand{\leq}{\leqslant}
\renewcommand{\geq}{\geqslant}
\newcommand{\noin}{\noindent}

\newcommand{\ep}{\varepsilon}
\newcommand{\ept}{\tilde{\ep}}
\renewcommand{\phi}{\varphi}

\newcommand{\mb}{\bm m }

\newcommand{\bk}{\bm k}
\newcommand{\bkt}{{\bk}}
\newcommand{\bkp}{\bm k}
\newcommand{\bkh}{\hat{\bk}}
\newcommand{\n}{{\bm n}}
\newcommand{\de}{\partial}
\newcommand{\x}{{\bm x}}
\newcommand{\bxi}{{\bm \xi}}

\renewcommand{\emptyset}{\varnothing}

\newcommand{\hpsi}{{\psi^\prime}}
\newcommand{\rf}[1]{(\ref{#1})}
\renewcommand{\rho}{\varrho}

\renewcommand{\kappa}{\varkappa}
\renewcommand{\setminus}{\smallsetminus}
\newcommand{\what}{\widehat}
\newcommand{\nonu}{\nonumber}
\newcommand{\oh}{\frac{1}{2}}

\newcommand{\ta}{\tilde{a}}
\newcommand{\tdelta}{\tilde{\delta}}
\newcommand{\z}{{\bm 0}}

\newcommand{\Am}{\bm A}

\renewcommand{\L}{\bm L}

\newcommand{\Ho}{H^\oh}
\newcommand{\Ht}{H^\frac{3}{2}}

\newcommand{\M}{\bm S}
\newcommand{\Mt}{\widetilde{\bm S}}
\newcommand{\Mi}{\bm M}

\newcommand{\Om}{\Omega}
\newcommand{\Omh}{\what{\Omega}}

\newcommand{\R}{{\mathbb R}}
\newcommand{\Z}{{\mathbb Z}}

\newcommand{\Id}{\bm I}

\renewcommand{\Im}{\bm I}
\newcommand{\J}{\bm J}

\renewcommand{\O}{{\cal O}}

\newcommand{\N}{\mathsfbfit N}
\newcommand{\Np}{{\N}^{+}}
\newcommand{\Nm}{{\N}^{-}}

\newcommand{\gi}{\gamma_{-}}
\renewcommand{\ge}{\gamma_{+}}
\newcommand{\om}{\omega}

\newcommand{\ei}{\varepsilon}

\newtheorem{theorem}{Theorem}[section]

\newtheorem{lemma}{Lemma}[section]

\newtheorem{remark}{Remark}

\addtocounter{figure}{0}
\newcommand{\I}{\ensuremath{\mathrm{i}}}
\newcommand{\E}{\ensuremath{\mathrm{e}}}
\newcommand{\D}{\ensuremath{\mathrm{d}}}

\begin{document}

\title{Local gaps in three-dimensional periodic media}
\author{
Yuri~A. Godin\thanks{Email: ygodin@uncc.edu} and Boris Vainberg\thanks{Email: brvainbe@uncc.edu} \\
The University of North Carolina at Charlotte, \\
Charlotte, NC 28223 USA
}

\maketitle

\begin{abstract}
We consider the propagation of acoustic waves in a medium with a periodic array of small inclusions of arbitrary
shape. The inclusion size $a$ is much smaller than the array period. 
We show that global gaps do not exist if $a$ is small enough. The notion of local gaps which depends on the choice of the
wave vector $\bk$ is introduced and studied. We determine analytically the location of local gaps for the Dirichlet
and transmission problems.
\end{abstract}



\bigskip

\section{Introduction}
\setcounter{equation}{0}

The discovery of photonic and phononic crystals, the creation and advancement of electromagnetic and acoustic metamaterials stimulated the development of efficient numerical methods for calculating their bandgap structure. The main tool for the evaluation of dispersion relations was the plane wave expansion method \cite{Kushwaha:94, Johnson:2001, Wu:2004}. Later, other methods such as the finite-difference
time-domain \cite{Fan:1996, Tanaka:2000, Taflove:2000}, the finite element \cite{Kuchment:99, Dobson:1999, Giani:2012} and boundary element
methods \cite{JJWM:11,Li:2013} employed. All these methods reduce the bandgap problem to the solution of an eigenvalue problem from which propagating frequency is determined.
Analytic determination of bandgaps is a more challenging problem which becomes possible by introducing small or large parameters. This situation occurs in solid-state physics in the study of the energy bands of electrons when
the periodic potential of crystals is small \cite{AM:76, Kittel:04}. Other examples are dealing with
high contrast materials for inclusions and the matrix, such results were established in \cite{Figotin:1996b, Ammari:2009, Lipton:2017, Lipton:2022a}.

The aim of standard bandgap analysis is
the determination of absolute gaps consisting of frequencies for which no waves propagate in any direction.
We will show the absence of absolute gaps in our problem and investigate {\it local gaps} defined by pairs $(\om,\bk_0)$ for which waves with frequency $\om$ do not propagate for small variations of the Bloch vector $\bk_0$. A detailed definition of a local gap is given later.

We study the propagation of Bloch 
waves in a medium containing a periodic array of identical inclusions. The amplitude $u$ of the waves is governed by the
equation
\begin{align}
 \Delta u + k^2 u =0,
 \label{u}
\end{align}
where $k = \om/c$ is the wave number, $\om$ is the time frequency, and $c$ is the wave propagation speed (different in the host medium and inclusions). Solution of \rf{u} is sought in the form
\begin{align*}
 u(\x) = \Phi (\x) \E^{-\I \bkp \cdot \x},
\end{align*}
where $\bk = (k_1, k_2, k_3)$ is the Bloch vector, and function $\Phi(\x)$ is periodic with the period of the lattice. The representation above is equivalent to
\begin{align}
 \rrbracket  \E^{\I\bkp \cdot \x} u(\x) \llbracket =0,
 \label{B1}
\end{align}
where $\rrbracket f \llbracket$ denotes the jump of $f$ and its gradient across the opposite sides of the cells of periodicity.

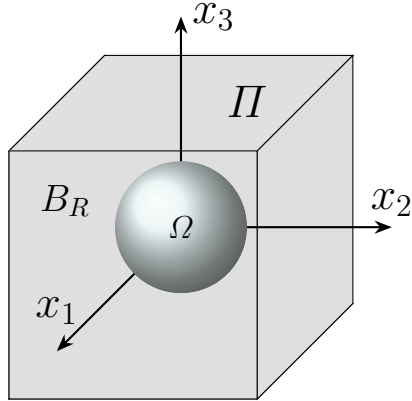
\begin{figure}[th]
\begin{center}
\begin{tikzpicture}[scale=1.1,>=Stealth]
\def\l{1.5}
  \draw[fill=lightgray!50] (-\l,-\l,-\l) -- (-\l,-\l,\l) -- (\l,-\l,\l) -- (\l,-\l,-\l) -- cycle; 
   \draw[fill=lightgray!50] (-\l,-\l,\l) -- (-\l,-\l,-\l) -- (-\l,\l,-\l) -- (-\l,\l,\l) -- cycle; 
   \draw[fill=lightgray!50] (-\l,\l,-\l) -- (-\l,\l,\l) -- (\l,\l,\l) -- (\l,\l,-\l) -- cycle; 
   \draw[fill=lightgray!50] (-\l,-\l,\l) -- (\l,-\l,\l) -- (\l,\l,\l) -- (-\l,\l,\l) -- cycle; 
   \draw[fill=lightgray!50] (\l,-\l,\l) -- (\l,-\l,-\l) -- (\l,\l,-\l) -- (\l,\l,\l) -- cycle; 

\draw [->, thick] (0,0,0) -- (0,0,2.6*\l);
\node [above] at (0,0.2,2.6*\l) {\Large  $x_1$};
\draw [->, thick] (0,0,0) -- (0,1.7*\l,0) node [right] {\Large  $x_3$};
\draw [->, thick] (0,0,0) -- (1.7*\l,0,0) node [above] {\Large  $x_2$};
 \shade[ball color = teal!10!white, opacity = 0.9] (0,0,0) circle (0.8cm);

 \def\eggheight{3mm}
 \begin{turn}{45}
  \path[ball color=white, opacity = 1]
  plot[domain=-pi:pi,samples=100,shift={(0,0)}]
  ({.8*\eggheight *cos(\x/4 r)*sin(\x r)},{-\eggheight*(cos(\x r))})
  -- cycle;
  \end{turn}
  \node at (0,0) {$\Om$};
  \node [above] at (0.8,1.2,0) {\Large$\Pi$};
  \node [above] at (-1.4,0,0) {\large$B_R$};
\end{tikzpicture}
\end{center}
\caption{
The cell of periodicity $\Pi$ with a small inclusion $\Om$ of size $a$. A ball $B_R \subset \Pi$ of radius $R$ centered at the origin encloses the inclusion $\Om$.
}
\label{cell}
\end{figure}

For simplicity, we assume that the fundamental cell of periodicity $\Pi$ is a cube $[-\pi,\pi]^3$. Denote by $\Om$ the domain occupied by the inclusion in $\Pi$ (see Figure \ref{cell}).
We select the origin for the coordinate system to be in $\Om$ and suppose that $\Om$ has a small size. Specifically, $\Om = \Om (a)$ is produced by compressing an $a$-independent region $\what{\Om}$ of arbitrary shape using the factor $a^{-1}$, $0 < a \ll 1$. In other words, the transformation $\x \to a\bxi $ maps  $\Om (a) \subset \R^3_{\x}$ into $\Omh \subset \R^3_{\bxi}$. We assume that the boundary $ \de\Om$ belongs to the class  $C^{1,\beta}$, which means that the functions that define the boundary have first-order derivatives that belong to the H\"{o}lder space with index $\beta$.

We consider first the case when homogeneous Dirichlet conditions, which arise for some soft inclusions, are imposed on the boundary. The transmission problem is considered in \S \ref{trans}. Then function $u$ in the fundamental cell
of periodicity $\Pi$ obeys the equation
\begin{align}
\label{Hz1}
 \Delta u + k^2 u &=0, \quad u \in H^2 (\Pi \setminus \Om)
\end{align}
and the boundary conditions
\begin{align}
  \left. u \right|_{\de\Om} &= 0, \quad  \rrbracket  \E^{\I\bk \cdot \x} u(\x) \llbracket =0.
 \label{bc1a}
\end{align}

A non-trivial solution of \rf{Hz1},\rf{bc1a} with fixed $\Omh$ exists not for all values of $k=\om/c, a,$ and $\bk$. The relation between the parameters $\om, a, \bk$, for which there is a non-trivial solution, is called the dispersion relation. The problem \rf{Hz1},\rf{bc1a} is one-periodic in each component of $\bk$, and the dispersion relation is a periodic function defined in $\R^3_{\bk}$.

The goal of this work is to find gaps i.e. the values of time frequency $\om=k c$ such that \rf{Hz1},\rf{bc1a} has only trivial solutions for all Bloch vectors $\bk\in\R^3$. Thus, the waves with frequency $\om$ do not propagate if $\om$ belongs to a gap. These gaps are also called global or complete gaps.
Geometrically, global gaps are empty intervals in the orthogonal projection of the dispersion surface onto the vertical $\om$-axis. Global gaps in $\om$ correspond to the gaps in the spectrum $\lambda=k^2$ of $-\Delta$ in $\Pi\setminus\Om(a)$ with an arbitrary $\bk$ in the boundary condition \rf{bc1a}.
We will also consider {\it local gaps} whose rigorous definition will be provided later in this section.

\begin{figure}[ht]
\centering
  \includegraphics[width=.6\linewidth]{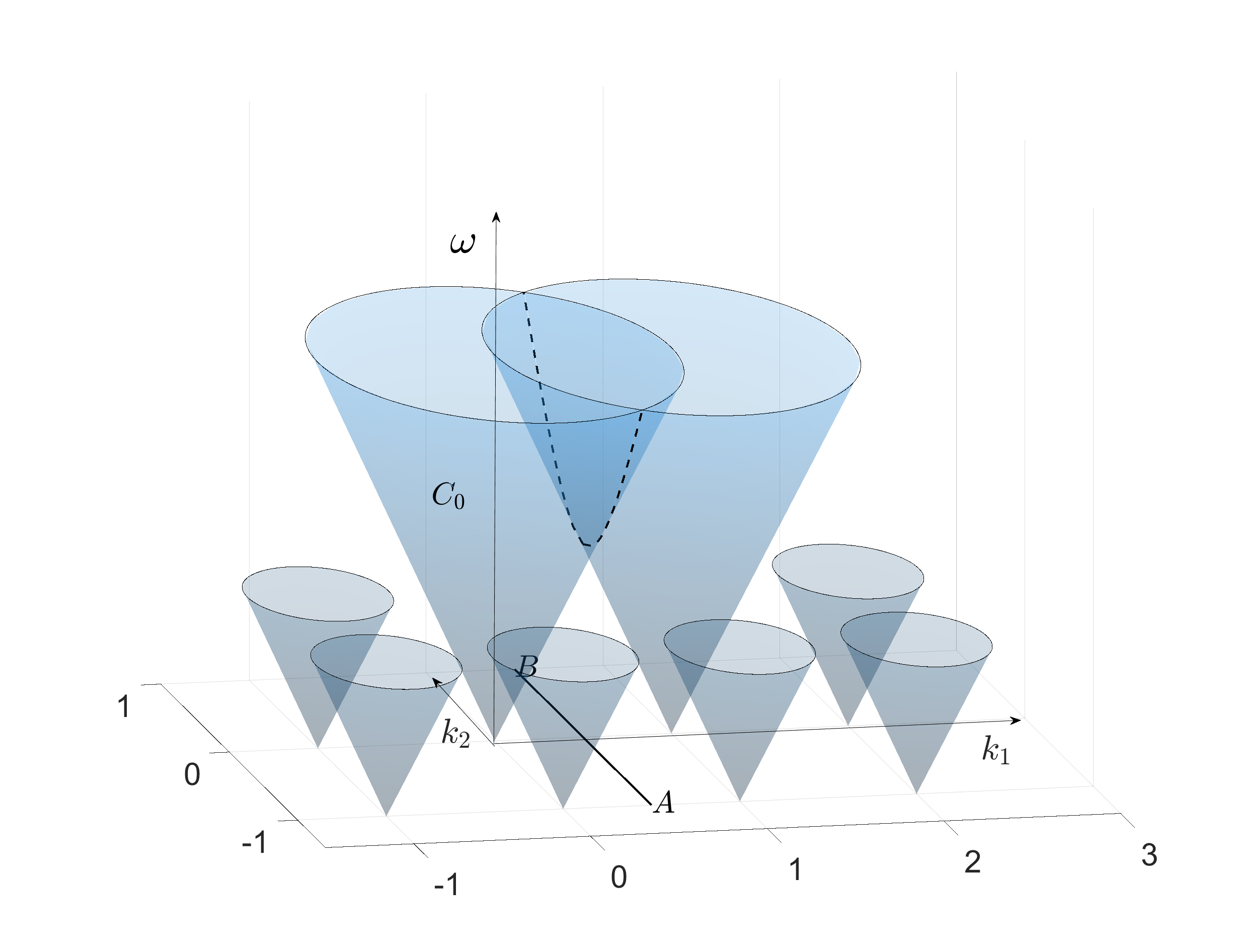}
 \caption{Several cones forming the dispersion surface of the unperturbed problem in dimension two.}
  \label{cones}
\end{figure}

Consider the unperturbed problem.
Observe that in the absence of inclusions, function $\E^{-\I \bk \cdot \x}$ satisfies \rf{Hz1},\rf{bc1a}  with $k = |\bk|$, i.e.
$\om = c|\bk|$.
For any vector $\mb = (m_1 , m_2,  m_3 )$ with integer components, functions $ \E^{-\I(\bkp - \mb) \cdot \x}$ are also satisfy the same boundary condition \rf{bc1a} and equation  \rf{Hz1} with $\om =c |\bk - \mb|$.
The graph of the infinite-valued periodic in $\bk$ dispersion function for the unperturbed problem consists of the set of all the cones $\om = c|\bk - \mb|$, see Figure \ref{cones}.
\vspace{-10mm}

For some values of $\bk$, the unperturbed problem \rf{Hz1},\rf{bc1a} with $k = |\bk|$ has solutions different from $ \E^{-\I\bk \cdot \x}$, and the solution space has dimension $n >1$. All such solutions have the form $ \E^{-\I(\bk - \mb) \cdot \x}$ provided that $|\bk| = |\bk - \mb|$ and $\mb$ is a vector with integer components.
The values of $\bk$ for which $n > 1$ are called {\it exceptional of order $n$}.
 Geometrically, $n$ refers to the number of points $\mb$ in the integer lattice $\mathbb Z^3$ (including the origin) belonging to the Ewald sphere \cite{Kittel:04} centered at $\bkp$ with radius $|\bkp|$, i.e. equation
 \begin{align}
 |\bk| = |\bk - \mb|
 \label{Ew}
 \end{align}
 has $n>1$ solutions $\mb \neq \z$ for exceptional $\bk$, and only one solution $\mb = \z$ for non-exceptional $\bk$.
Equation \rf{Ew} can be written as
\begin{align}
2\bk \cdot \mb = |\mb|^2, \quad \mb \in {\mathbb Z}^3 \setminus \z,
\label{Ew1}
\end{align}
i.e. exceptional points belong to planes in ${\mathbb R}^3_{\bk}$ whose distance from the origin goes to infinity as $|\mb| \to \infty$. An exceptional point $\bk$ has order $n > 1$ if it belongs to $n-1$ planes \rf{Ew1}. One can also describe
exceptional wave vectors as values of $\bk$ for which $n>1$ cones in Figure \ref{cones} interest each other.
In particular, in the two-dimensional analog of our problem in Figure \ref{cones},
the line $AB$  consists of exceptional points of order at least two.

The set of solutions of the problem  \rf{Hz1},\rf{bc1a} remains the same if the wave vector is changed by a period of the reciprocal lattice.
Thus it is enough to study local gaps related to one cone centered at the origin: $C_0 :=\{\om=c|\bk|\}$. This will allow us to define local gaps related to other cones in Figure \ref{cones} using the periodicity of solutions in $\bk$.

Let us fix a wave vector $\bk_0$.
{\it The local gap} $g(\bk_0)$ corresponding to $\bk_0$ and associated with $C_0$ consists of frequencies $\om = \om_a$ for which
Bloch waves do not propagate when wave vector $\bk$
has the same direction as $\bk_0$ and $|\bk - \bk_0|$ is sufficiently small. More precisely, Bloch
waves with the frequency $\om_a$ do not exist
when $\bk = (1+\delta)\bk_0$ for sufficiently small $a$-independent $|\delta|$.

It is not expedient to consider frequencies $\om$
 if the point $(\om,\bk_0)$ is located at a positive, $a$-independent distance $d$ from the set of cones shown in Fig. 2 (dispersion surfaces of the unperturbed problem). Bloch waves with parameters $(\om,(1 + \delta)\bk_0)$ and $\delta < d/2$
do not propagate in the unperturbed problem and the problem with inclusions if $a$ is small enough.
Thus, we introduce one more requirement in the definition of local gap $g(\bk_0)$. We assume that $\om=\om_a \to \om_0$ as $a \to 0$ where $(\om_0, \bk_0) \in C_0$, i.e. $\om_0 = c|\bk_0|$.

Let $\om_a$ tends to $\om_0= c|\bk_0|$. We will show that
\begin{enumerate}
\item[(a)]
If $(\om_0, \bk_0)$ belongs to a cone $C_0$ but not to the intersection of the cones, then $g(\bk_0) = \emptyset$,
i. e. there are no local gap for wave vector $\bk_0$ and frequencies $\om$ near $\om_0$.

 \item[(b)]
 If $(\om_0, \bk_0)$ belongs to an intersection of $C_0$ with some other cone shown in Figure \ref{cones} then local gap $g(\bk_0)$ exists in a neighborhood
 of $\om_0$ for small $a$ if and only if $\dst \frac{|\bk_0|}{|\mb_0|} < \frac{\sqrt{2}}{2}$, where $\bk_0, \mb_0$ satisfy \rf{Ew}. The location of the gap will be specified. We do not consider the points $(\om_0, \bk_0)$ that belong to the intersection of three or more cones and those for which $\dst \frac{|\bk_0|}{|\mb_0|} = \frac{\sqrt{2}}{2}$.

 \item[(c)]
 We will show that global gaps do not exist in any fixed interval $\epsilon<\om<\epsilon^{-1}$ of the time frequency $\om$ if the size $a$ of the inclusion is small enough.
\end{enumerate}

 Note that the set of values of $\bk$ omitted from the consideration in item (b) forms a one-dimensional manifold in $\R^3$. Indeed, the set of exceptional vectors of multiplicity $n\geq 2$ is given by the planes \rf{Ew1}. Omitted from the consideration set of exceptional vectors of order higher than two are lines of intersections of those planes. Also, condition $\dst \frac{|\bk_0|}{|\mb_0|} = \frac{\sqrt{2}}{2}$ defines curves on planes \rf{Ew1}, see Figure 3, where these curves are circles that are boundaries of shaded disks.

\section{Outline of the approach}
\setcounter{equation}{0}
\label{out}

In \cite{GV:22, GV:23} we devised a new approach to determining the dispersion relation which allowed us to find the asymptotic expansion of the dispersion relation as $a\to 0$. The expansions have a form of power series in $a$ which are different for non-exceptional and exceptional wave vectors $\bk$. The expansion depends analytically on $\bk$ when $\bk$ is not exceptional, see a discussion below. To study gaps, we need to develop further these results in the present paper to obtain asymptotics in $a$ which is uniform in $\bk$ in a neighborhood of exceptional points.

As before, we enclose the inclusion $\Om$ in a ball $B_R \subset \Pi$ of radius $R > a$ centered at the origin and split $\Pi$ into the ball $B_R$ and its complement $\Pi \setminus B_R$ (see Figure \ref{cell}).
Next, we consider the following two problems in these two regions:
\begin{align}
 \left( \Delta + k^2 \right) u(\x) &= 0, \quad \x \in (B_R \setminus \Om), \quad  \left. u \right|_{\de\Om}= 0, \quad
 \left. u\right|_{r=R}=\psi, \label{uin} \\[2mm]
 \left( \Delta + k^2 \right) v(\x) &= 0, \quad  \x\in \Pi \setminus B_R,\quad  \left\rrbracket  \E^{\I \bkp \cdot \x} v(\x) \right\llbracket =0, \quad \left. v\right|_{r=R}=\psi,  \label{uout}
\end{align}
where
\begin{align}
 k &= \om/c = (1+\ep) |\bk|,
 \label{k}
\end{align}
$\ep$ will be determined below, $u \in H^2 (B_R \setminus \Om), ~ v \in H^2 (\Pi \setminus B_R)$.
The unperturbed problems with $a=\ep=0$ are uniquely solvable for smooth enough $\psi$ for all $R$ except a countable set $\{R_i\}$,
 and we choose $R \notin \{R_i$\}.  Then solutions of \rf{uin},\rf{uout} defines two Dirichlet-to-Neumann operators
$\Nm_{a,\ep}$, $\Np_{\bkp,\ep}$: $H^{3/2} (\de B_R) \to H^{1/2} (\de B_R)$, i.e.
\begin{equation}
\label{}
\Nm_{a,\ep}  \psi = \left.\frac{\de u}{\de r}\right|_{r=R}, \quad  \quad
\Np_{\bkp,\ep}  \psi = \left.\frac{\de v}{\de r}\right|_{r=R}.
\end{equation}
The existence of a non-trivial solution to the problem \rf{Hz1},\rf{bc1a} is equivalent to the existence of a zero eigenvalue (EZE) for the operator $\Np_{\bkp,\ep} - \Nm_{a,\ep}$. This abbreviation pertains specifically to the operator $\Np_{\bkp,\ep} - \Nm_{a,\ep}$. Hence, the dispersion relation is given by \rf{k} with $\ep=\ep (a,\bk)$ determined from the EZE condition. In \cite{GV:22, GV:23}, the reduction of the dispersion relation to the EZE problem enabled us to use regular perturbation techniques to study dispersion relations for each $\bk$.

Observe that the choice of $R$ and therefore the definition of operators $\N^\pm$ depend on $\bk$.
Since $k^2$ is an analytic function of $\bk$ and since $R$ is chosen in such a way that the problem  \rf{uout} is uniquely solvable when $\bk = \bk^\prime$ is fixed and $\ep=0$,
the operator $\Np_{\bkp,\ep}$ is analytic in $\bk$ and $\ep$ when $|\bk - \bk^\prime| + |\ep| \ll 1$. We will call this property the local analyticity in $\bk \in \R^3$ and small $\ep$.
Similarly, the operator $\Nm_{a,\ep}$ is locally analytic in $\bk$ and small $\ep$.
The smoothness of the operator $\Nm_{a,\ep}$ in $a$ is not obvious, because
problem \rf{uin} is a singular perturbation of the corresponding problem without the inclusion ($\Om = \emptyset$), and solutions of \rf{uin} do not have uniform limits as $a \to 0$. However, the solution of \rf{uin} is infinitely smooth in $a$ when $a$ is small and $\x$ is outside of a fixed neighborhood of the origin, and it was shown in \cite{GV:22, GV:23} that the interior DtN operator
$\Nm_{a,\ep}$ is infinitely differentiable in $a,~a\ll1$.
We summarize these properties in the Lemma.
\begin{lemma}
\label{l21}
 \begin{enumerate}
  \item
  The operator $\Np_{\bkp,\ep}$ is locally analytic in $\bk \in \R^3$ and $\ep$.
  The operator $\Nm_{a,\ep}$ is locally analytic in $\bk$ and $\ep$ and infinitely smooth in $a,~a\ll1$.
  \item
  The dispersion relation \rf{k} for which  \rf{Hz1},\rf{bc1a} have a nontrivial solution is given by the same function $\ep = \ep(a,\bk)$
  for which the EZE condition holds.
 \end{enumerate}
\end{lemma}

The Lemma reduces the study of the dispersion relation to the EZE problem for a smooth operator function.
Special matrix representations of operators $\Np_{\bkp,\ep},~  \Nm_{a,\ep}$ were constructed in \cite{GV:22, GV:23} that allowed us to solve the EZE problem and find the asymptotics of the dispersion relation $\ep=\ep (a,\bk)$ as $a\to0$ which is analytic in $\bk$ for non-exceptional $\bk$ and is valid for fixed exceptional $\bk$. The gaps are related to exceptional points, and the study of gaps requires a uniform in $\bk$ asymptotics of $\ep (a,\bk)$. The main goal of the present paper is to construct a matrix representation of the operator $\Np_{\bkp,\ep}-  \Nm_{a,\ep}$, allowing to obtain uniform in $\bk$ asymptotics of the dispersion relations near the exceptional points and use them to study the gaps.

The study of the operator $\Np_{\bkp,\ep} - \Nm_{a,\ep}$ is simplified if we subtract $\Nm_{0,\ep}$ from both terms of the difference above.

\section{Necessary prerequisites}
\label{np}
\setcounter{equation}{0}

This section provides the results from \cite{GV:22, GV:23} needed to understand the goal and approach of the present paper.

The dispersion surface for the inclusionless ($\Omega=\emptyset, k = |\bk|$) problem \rf{Hz1}, \rf{bc1a} consists of the cone $C_0:~\om=c|\bk|$ and its shifts shown in Fig. 1,2. Since the solutions of \rf{Hz1}, \rf{bc1a} and of the corresponding transmission problem are periodic in $\bk$ and their dispersion surfaces are close to the surface for the inclusionless problem, it is enough to study the dispersion surface of the problems with inclusions only in a small neighborhood of $C_0$. The information about other parts of the surface can be obtained simply by the shift in $\bk$. Thus, similar to the approach used in  \cite{GV:22, GV:23}, we will assume first that $(\om,\bk)$ is located in a small neighborhood of $C_0$. \\
{\bf Non-exceptional wave vectors.}
Let $\bk$ be a non-exceptional wave vector. Then solution space of the unperturbed problem \rf{Hz1},\rf{bc1a} ($\Om = \emptyset, k = |\bk|$) is one-dimensional and consists of functions proportional to $\E^{-\I \bk \cdot \x}, ~\x \in \Pi$. Respectively, the kernel of the operator $\Np_{\bkp,0} - \Nm_{0,\ep}$ consists of functions proportional to
$\hpsi = \E^{-\I \bk \cdot \x}, ~\x \in \de B_R$ which is the restriction of the exponent on $\de B_R$. We present the domain $\Ht(\de B_R)$ and the range $\Ho (\de B_R)$ of the operator
$\Np_{\bkp,\ep} - \Nm_{a,\ep}$ as a direct sum of the one-dimensional space $\mathscr E$ of functions proportional to $\hpsi$ and spaces $\mathscr E_{\bot,d}$,$ \mathscr E_{\bot,r}$, respectively, of functions orthogonal in $L^2 (\de \Om)$ to $\hpsi$. Here subindexes ``$d$'' and ``$r$'' stay for the domain and the range.

We represent each element $\psi$ in the domain and the range of operator $\Np_{\bkp,\ep} - \Nm_{a,\ep}$ in the vector form $\psi=( \psi_{\mathscr E},\psi_\bot),$ where $\psi_{\mathscr E}$ is the projection in $L^2(\de B_R)$ of the function $\psi$ on $\hpsi$, and $\psi_\bot$ is orthogonal to $\hpsi$ in $L^2(\de B_R)$, and  use
matrix representation of all the operators in the basis $( \psi_{\mathscr E},\psi_\bot).$ Due to Green's formula, the operators $\Np_{\bkp,\ep},$ $\Nm_{a,\ep}$ are symmetric and therefore
\begin{equation}
\label{A}
 \Np_{\bk,0} - \Nm_{0,\ep} = \left(
 \begin{array}{cc}
  0 & 0 \\[2mm]
  0 & \Am
 \end{array}
 \right),
\end{equation}
where $\Am:\mathscr E_{\bot,d}\to \mathscr E_{\bot,r}$ is an isomorphism.

Since operator $\Np_{\bk,\ep}-\Nm_{0,\ep}$ depends  smoothly on $\ep$,
we have
\begin{equation}
\label{conc}
 \Np_{\bk,\ep} - \Nm_{0,\ep} = \left(
 \begin{array}{cc}
  C\ep +\O(\ep^2) &\O(\ep) \\[2mm]
  \O(\ei) &\Am +\O(\ep)
 \end{array}
 \right),
\end{equation}
where $C$ is a constant, and $\ep$ is defined in \rf{k}. In \cite{GV:22} it was shown that $C=2|\bk|^2 |\Pi|$.

Matrix representation for the interior operator has the form
\begin{equation}
\label{mm}
 \Nm_{a,\ep} - \Nm_{0,\ep} = \left(
 \begin{array}{cc}
  4\pi q a +\O(a^2 +a |\ep|) &\O(a + |\ep|) \\[2mm]
 \O(a + |\ep|) &O(a + |\ep|)
 \end{array}
 \right),
\end{equation}
where $qa$ depends on the shape of $\de \Om$ and is equal to its electrical capacitance, i.e. total charge of $\de \Om$ when its potential is unity. In particular, $q=1$ for a sphere. Thus
\begin{equation}
\label{conc1}
 \Np_{\bk,\varepsilon} - \Nm_{a,\varepsilon} = \left(
 \begin{array}{cc}
  C\ep -4\pi q a +\O(a^2 +\ep^2) &\O(a+|\ep|) \\[2mm]
 \O(a+|\ep|) &\Am +\O(a+|\ep|)
 \end{array}
 \right), \quad C=2|\bk|^2 |\Pi|.
\end{equation}
This matrix representation allows one to find $\ep=\ep(a,\bk)$ for which the kernel of the operator $\Np_{\bk,\varepsilon} - \Nm_{a,\ep}$ is not empty. Using \rf{conc1}, we write the equation  $(\Np_{\bk,\varepsilon} - \Nm_{a,\varepsilon} )\psi=0$ in the form
\begin{equation}
\label{conc5}
  \left(
 \begin{array}{cc}
  2\ep |\bk|^2 |\Pi| -4\pi q a +\O(a^2 +\ep^2) &\O(a+|\ep|) \\[2mm]
 \O(a+|\ep|) &\Am +\O(a+|\ep|)
 \end{array}
 \right)
 \left(
 \begin{array}{cc}
 \psi_{\mathscr E } \\[2mm]
 \psi_\bot
 \end{array}
 \right)=0.
\end{equation}
Due to the invertibility of $\Am$, we solve equation above for $\psi_\perp:~ \psi_\bot=\O(a+|\ep|) \psi_{\mathscr E}$, and reduce \rf{conc5} to a simple equation for $\psi _{\mathscr E}$. This implies that function $\ep(a,\bk)$ can be found by equating the top left element in \rf{conc5} with a different remainder term to zero. Thus,  the implicit function theorem implies the following statement:
\begin{lemma}\label{l22}
Let $\bk'\neq \z$ be an arbitrary non-exceptional wave vector and $k'=|\bk'|$. Then
 \begin{enumerate}
  \item
 The dispersion relation $\om=c(1+\ep)|\bk|$ is uniquely defined in a neighborhood of the point $(k',\bk')\in \R^4$ if $a$  is sufficiently small. Function $\ep = \ep(a,\bk)$ is  infinitely smooth in $a, \bk$ and analytic in $\bk$ when $\bk \neq \z$.
  \item
The following asymptotics holds
\begin{equation}
\label{asep}
 \ei=\frac{2\pi qa}{|\bk|^2|\Pi|}+O(a^2),\quad a\to 0, \quad \bk \neq \z.
\end{equation}
 \end{enumerate}
\end{lemma}
\noin
{\bf Exceptional wave vectors of order two.}
Let the Bloch vector $\bk=\bk_0$ in \rf{bc1a} be exceptional of multiplicity two, and therefore there exists a non-zero vector $\mb$ with integer components such that $|\bk_0|=|\bk_1|$ for $\bk_1=\bk_0-\mb$. Then solution space of the unperturbed problem \rf{Hz1},\rf{bc1a} is spanned by functions $\E^{-\I \bk _0\cdot \x}, ~ ~\E^{-\I \bk_1 \cdot \x}, ~\x \in \Pi$. Respectively, the kernel $\mathscr E$ of the operator $\Np_{\bkp,0} - \Nm_{0,\ep}$ is two-dimensional and consists of linear combinations of functions
$\psi_i^\prime = \E^{-\I \bk_i \cdot \x}, ~\x \in \de B_R, ~i=0,1$. Once again we present the domain and the range of the operator
$\Np_{\bkp,\ep} - \Nm_{a,\ep}$ as a direct sum of space $\mathscr E$ and its orthogonal in $L^2 (\de \Om)$ complements, and write all the operators in corresponding matrix form.

The formula \rf{A} remains valid when $\bk = \bk_0$ is an exceptional wave vector of order two if zeroes of the matrix on the right-hand side of \rf{A} are understood as the zero elements of the first two rows and first two columns of the matrix. In addition, formulas \rf{conc} and \rf{mm} remain valid if constants $C$ and $q$ are replaced by $C\Id$ and $q\J$, respectively, where $\Id$ and $\J$
are $2\times2$ identity and all-ones matrices, respectively.
The analogue of \rf{conc1} has the following $2\times2$ block $\M$ in the top left corner of the matrix:
\begin{align}
\M=2\ep |\bk_0|^2 |\Pi|\Id - 4\pi q a \J  + \O(a^2 + \ep^2),
\label{M1}
\end{align}
and therefore, \rf{conc5} now takes the form
\begin{equation}
\label{conc6}
  \left(
 \begin{array}{cc}
 \M +\O(a^2 +\ep^2) &\O(a+|\ep|) \\[2mm]
 \O(a+|\ep|) &\Am +\O(a+|\ep|)
 \end{array}
 \right)
 \left(
 \begin{array}{cc}
 \psi_{\mathscr E } \\[2mm]
 \psi_\bot
 \end{array}
 \right)=0.
\end{equation}
The arguments after \rf{conc5} imply
that the dispersion relation at $\bk=\bk_0$ is defined by the equation $\det{\Mt}=0$ where ${\Mt}$ and $\M$ have the same asymptotics with different specific values of the remainders. Since the eigenvalues of $\J$ are distinct, the implicit function theorem can be applied to the equation $\det{\Mt}=0$ for the variable $\ep/a$. This leads to the existence of the following two roots:
\begin{equation}
\label{expp}
\ep=\ep_1=\frac{2\pi qa}{|\bk_0|^2|\Pi|}+O(a^2),  \quad   \ep=\ep_2=O(a^2),~~ a\to0.
\end{equation}

In addition, the following important observation was made in \cite{GV:22}: there are no solutions to the problem \rf{Hz1},\rf{bc1a} with $\bk=\bk_0$ propagating in one direction of either vector $\bk_0$ or $\bk_1$. Each of the Bloch solutions with $\ep=\ep_i, ~i=0,1,$ is a cluster of waves propagating in both directions.

\section{Asymptotics of the DtN operators in the neighborhoods of exceptional points}
\setcounter{equation}{0}

The results concerning non-exceptional points described above are valid for {\it all} non-exceptional points. In particular, the function \rf{asep} is analytic in $\bk$. On the contrary, formulas \rf{expp} were obtained for {\it fixed} exceptional points whose neighborhoods contained mostly non-exceptional points. We need to extend \rf{A}-\rf{conc1} in such a way that the corresponding asymptotics would be valid in the neighborhood of fixed wave vector $\bk=\bk_0$ of multiplicity two.

We do not consider the whole neighborhood of $\bk_0$. It will be sufficient to focus only on the values of $\bk$ that belong to
the ray through $\bk_0$. Thus $\bk = (1+\delta) \bk_0$, ~$|\delta| \ll 1$.
We will denote this segment by  $I$. It contains one exceptional point $\bk=\bk_0$, and other points in $I\setminus \bk_0$ are non-exceptional if $\delta$ is small enough. Matrix representation \rf{A}-\rf{conc1} in the previous section was based on one-dimensional space $\mathscr E$ when $\bk\in I\setminus \bk_0$ and two-dimensional space $\mathscr E$ when $\bk=\bk_0$.
Our next goal is to obtain the dispersion relation for the entire interval $I$ using the two-dimension space
$\mathscr E$ spanned by functions $\psi_0$ and $\psi_1$.
This will allow us to obtain the asymptotics of the dispersion relation as $a \to 0$, $\bk \in I$,  which is uniform in $\bk$.

We start with a matrix representation of $\Np_{\bk,\ep} - \Nm_{0,\ep}$  for $\bk \in I$ similar to \rf{A}-\rf{conc}. As mentioned at the end of the previous section, if $\ep=\delta=0$ then formula \rf{A} remains valid for the two-dimensional space $\mathscr E$.
Since all the operators involved depend infinitely smoothly on small $\ep, \delta$ and the space $\mathscr E$ is independent on
$\ep, \delta,$
the matrix representation of $\Np_{\bk,\ep} - \Nm_{0,\ep}$ for sufficiently small $\ep$ and $\delta$ has the form
\begin{align}
\label{conc2}
 \Np_{\bk,\ep} - \Nm_{0,\ep} = \left(
 \begin{array}{cc}
  \L(\ep,\delta) + \O(\ep^2 + \delta^2)  & \O(|\ep| + |\delta|) \\[2mm]
  \O(|\ep| + |\delta|) &\Am + \O(|\ep| + |\delta|)
 \end{array}
 \right), \quad \bk \in I,
\end{align}
where $\L$ is a $2\!\times \!2$ matrix which is linear in $\ep$ and $\delta$, and the entries of the matrix on the right-hand side are infinitely smooth in $\ep$ and $\delta$. The next Lemma reveals the dependence of $\L$ on $\ep$ and $\delta$.
\begin{lemma}
\label{4.1}
Matrix $\L$ has the form
 \begin{align}
 \L = 2\ep |\bk_0|^2 |\Pi| \Id
 + \delta |\Pi| \left[
 \begin{array}{cc}
  0 & 0 \\[2mm]
  0 & |\mb_0|^2
 \end{array}
\right],
\label{L1}
\end{align}
where the couple $(\bk_0, \mb_0)$ satisfies \rf{Ew}, \rf{Ew1}.
\end{lemma}
\noin
We will prove the lemma after we prove the following theorem.
\begin{theorem}
\label{t41}
 The matrix representation of the operator $\Np_{\bk,\ep} - \Nm_{a,\ep}$
 has the form
 \begin{equation}
\label{conc3}
 \Np_{\bk,\varepsilon} - \Nm_{a,\varepsilon} = \left(
 \begin{array}{cc}
  \M  & \O(a+|\ep|+ |\delta|)\\[2mm]
 \O(a+|\ep|+ |\delta|) & ~~\Am +\O(a+|\ep|+ |\delta| )
 \end{array}
 \right), \quad \bk \in I,
\end{equation}
where the $2\times 2$ upper left block is given by
 \begin{align}
\M=2\ep |\bk_0|^2 |\Pi| \Im - 4\pi a q \J
+ \delta  |\Pi| \left[
 \begin{array}{cc}
  0 & 0 \\[2mm]
  0 & |\mb_0|^2
 \end{array}
\right]
+ \O(a^2 + \delta^2 + \ep^2 ), \quad a+|\ep|+|\delta|\ll1,
 \label{tau}
\end{align}
where the couple $(\bk_0, \mb_0)$ satisfies \rf{Ew}, \rf{Ew1}  and $q$ is defined in \rf{mm}.
\end{theorem}
\begin{proof}{}
To obtain \rf{conc3} we need to subtract the matrix representation
of the operator $\Nm := \Nm_{a,\ep} - \Nm_{0,\ep}$, $\bk \in I$, from \rf{conc2}. $\Nm$ depends infinitely smoothly on $a,\ep,$ and $\bk$ and
can be considered as an infinitely smooth operator-function of $a,\ep,$ and $\delta$. Obviously, the operator is zero when $a=0$, and therefore can be written as $\Nm=a\N'+O(a(a+|\ep|+|\delta|))$, where operator $\N'$ does not depend on $a,\ep,\delta$. This formula, \rf{conc2} and Lemma \ref{L1} will justify \rf{conc3} if we show that the upper left $2\times 2$ block of the matrix representation for $a\N'$ is equal to $4\pi a q \J $. This block of $a\N'$ coincides with the linear in $a$ part of a similar block of the matrix representation for $\Nm$ with $\ep=\delta=0$. The latter one was obtained in \cite{GV:23} when the matrix representation was studied for a fixed exceptional wave vector $\bk=\bk_0$ of multiplicity two, see \S \ref{out}, and it equals $4\pi a q \J $.
\end{proof}

\renewcommand*{\proofname}{Proof of Lemma \ref{conc2}}
\begin{proof}
Since $\L = \ep \L_1 + \delta \L_2$ and formula \rf{conc2} when $\delta=0$ coincides with the matrix representation constructed in the previous section when $\bk= \bk_0$, we have $\L_1 = 2 |\bk_0|^2 |\Pi| \Id$.
Thus we need to justify only the form of $\L_2$ in \rf{L1} and therefore we set $\ep = 0$ in the calculations below.

The entries of $\L_2 = \{L_{i,j}\}$ are given by
\begin{align}
 L_{i,j} = \int_{\de B_R} (\Np_{\bk,0} - \Nm_{0,\ep}) \E^{-\I \bk_i \cdot \x}\cdot \E^{\I \bk_j \cdot \x}\, \D S
  + \O(\delta^2), \quad 0\leq i,j \leq 1,
 \label{L}
\end{align}
where $\bk_0, \bk_1$ are two related exceptional wave vectors of multiplicity two, i.e., $\bk_1 = \bk_0 - \mb_0$ and $|\bk_0| = |\bk_1|$.

Consider the solution $u = u_i$ of the problem in $\Pi \setminus \de B_R$ which is a union of the exterior problem in $\Pi \setminus B_R$ and the problem in $B_R$ without the inclusion:
\begin{align}
 \label{k-egv}
 \left( \Delta + |\bk|^2 \right) u_i(\x) &= 0, \quad  \x\in \Pi, \quad \x \notin  \de B_R,~~  \left. u_i\right|_{r=R}=\psi_i = \E^{-\I \bk_i \cdot \x},
 \quad i=0,1,
\end{align}
with the Bloch boundary conditions on $\de \Pi$
\begin{align}
 \label{BC}
 \left\rrbracket  \E^{\I \bkt \cdot \x} u_i(\x) \right\llbracket &=0.
\end{align}
Here we replaced $k$ by $|\bk|$ since $\ep = 0$.
We write \rf{L} in terms of $u_i$:
\begin{align}
 L_{i,j} = \int_{\de B_R} \left(\left.\frac{\de u_i}{\de r}\right|_{r=R+0} - \left.\frac{\de u_i}{\de r}\right|_{r=R-0}\right) \E^{\I \bk_j \cdot \x}\,\D S  + \O(\delta^2).
 \label{Lij}
\end{align}
Similar to solutions of \rf{uin}, \rf{uout}, functions  $u_i$ are infinitely smooth function of $\delta$ and therefore
\begin{align}
 u_i(\x) = \E^{-\I \bk_i \cdot \x} + \delta \cdot v_i(\x) +\O(\delta^2).
 \label{uexp}
\end{align}
Then $v_i(\x)$ are solutions of the following problems:
\begin{align}
 \label{v}
 \left( \Delta + |\bk_0|^2 \right) v_i(\x) &= -2|\bk_0|^2 \E^{-\I \bk_i \cdot \x}, \quad  \x\in \Pi, \quad \x \notin  \de B_R,~~ \left. v_i\right|_{r=R}=0, \\[2mm]
 \label{bv1}
 \left\rrbracket  \E^{\I \bk_0 \cdot \x} v_i(\x) + \I \bk_0 \cdot \x \right\llbracket &=0, \quad i=0,1,
 \end{align}
and \rf{Lij} takes the form
\begin{align}
 L_{i,j} = \int_{\de B_R} \left(\left.\frac{\de v_i}{\de r}\right|_{r=R+0} - \left.\frac{\de v_i}{\de r}\right|_{r=R-0}\right) \E^{\I \bk_j \cdot \x}\,\D S + \O(\delta^2).
 \label{Lij1}
\end{align}

In what follows we omit subscript ``$i$'' from $v=v_i$ and add subscripts ``$+$'' and ``$-$'' to the restriction of $v$ to the domains $\Pi \setminus B_R$ and $B_R$, respectively.

To evaluate the integral $\dst \int_{\de B_R} \frac{\de v_{+}}{\de \n}\, \E^{\I \bk_j \cdot \x}\, \D S$, where $v$ is solution of \rf{v}-\rf{bv1} in $\Pi \setminus B_R$ we
apply Green's second identity to the functions $v_{+}(\x)$ and $ \E^{\I \bk_j \cdot \x}$ in the domain $\Pi \setminus B_R$, $j=0,1$:
\begin{align}
\int_{\Pi \setminus B_R}(\Delta v_{+} + |\bk_0|^2 v_{+})  \E^{\I \bk_j \cdot \x} \, \D V
&= \int_{\de \Pi} \left(\frac{\de v_{+}}{\de \n}\, \E^{\I \bk_j \cdot \x} -v_{+}\,\frac{\de \E^{\I \bk_j \cdot \x} }{\de \n}  \right) \D S
\nonu \\[2mm]
&- \int_{\de B_R} \left(\frac{\de v_{+}}{\de r}\, \E^{\I \bk_j \cdot \x} -v_{+}\,\frac{\de \E^{\I \bk_j \cdot \x} }{\de r}  \right) \D S.
 \label{G2}
\end{align}
From here we have
\begin{align}
 & \int_{\de B_R} \frac{\de v_{+}}{\de r}\, \E^{\I \bk_j \cdot \x}\, \D S =
 \int_{\de \Pi} \left(\frac{\de v_{+}}{\de \n}\, \E^{\I \bk_j \cdot \x} -v_{+}\,\frac{\de \E^{\I \bk_j \cdot \x} }{\de \n}  \right) \D S
 + 2|\bk_0|^2 \int_{\Pi \setminus B_R} \E^{\I (\bk_j -\bk_i) \cdot \x} \, \D V \nonu \\[2mm]
 &=\int_{\de \Pi} \left(\frac{\de v_{+}}{\de \n}\, \E^{\I \bk_j \cdot \x} -v_{+}\,\frac{\de \E^{\I \bk_j \cdot \x} }{\de \n}  \right) \D S
 + 2|\bk_0|^2 \int_{\Pi} \E^{\I (\bk_j -\bk_i) \cdot \x} \, \D V - 2|\bk_0|^2 \int_{B_R} \E^{\I (\bk_j -\bk_i) \cdot \x} \, \D V.
 \label{dvn}
\end{align}
Evaluation of the similar integral for $v_{-}$ yields
\begin{align}
 \int_{\de B_R} \frac{\de v_{-}}{\de r}\, \E^{\I \bk_j \cdot \x}\, \D S &=
 -2|\bk_0|^2 \int_{B_R} \E^{\I (\bk_j -\bk_i) \cdot \x} \, \D V.
 \label{dvi}
\end{align}
Thus,
\begin{align}
\int_{\de B_R} \left(\frac{\de v_{+}}{\de r} - \frac{\de v_{-}}{\de r}\right) \E^{\I \bk_j \cdot \x}\, \D S =
\int_{\de \Pi} \left(\frac{\de v_{+}}{\de \n}\, \E^{\I \bk_j \cdot \x} -v_{+}\,\frac{\de \E^{\I \bk_j \cdot \x} }{\de \n}  \right) \D S
 +2|\bk_0|^2 \int_{\Pi} \E^{\I (\bk_j -\bk_i) \cdot \x} \, \D V.
 \label{dvni}
\end{align}
In the middle integral, we make the substitution
\begin{align}
 v_{+}(\x) =  -\I \bk_0 \cdot \x \,\E^{-\I \bk_0 \cdot \x} + \tilde{v}(\x).
 \label{}
\end{align}
Then $\tilde{v}(\x)$ satisfies the homogeneous condition \rf{bv1}
and the middle integral in \rf{dvni} with  $\tilde{v}(\x)$ instead of $v_{+}$ vanishes, and \rf{dvni} becomes
\begin{align}
\int_{\de B_R} \left(\frac{\de v_{+}}{\de r} - \frac{\de v_{-}}{\de r}\right) \E^{\I \bk_j \cdot \x}\, \D S &= -
\int_{\de \Pi} \left(\frac{\de  (\I \bk_0 \cdot \x \,\E^{-\I \bk_i \cdot \x})}{\de \n}\, \E^{\I \bk_j \cdot \x} - \I \bk_0 \cdot \x \,\E^{-\I \bk_i \cdot \x}\,\frac{\de \E^{\I \bk_j \cdot \x} }{\de \n}  \right) \D S \nonu \\[2mm]
 &+2|\bk_0|^2 \int_{\Pi} \E^{\I (\bk_j -\bk_i) \cdot \x} \, \D V.
\end{align}
Using again Green's second identity we obtain
\begin{align}
 & \int_{\de B_R} \left(\frac{\de v_{+}}{\de r} - \frac{\de v_{-}}{\de r}\right) \E^{\I \bk_i \cdot \x}\, \D S \nonu \\[2mm]
 &=
 -\int_{\Pi} \left( \E^{\I \bk_j \cdot \x} \left(\Delta + |\bk_0|^2 \right) \left(\I \bk_0\cdot \x \,\E^{-\I \bk_i \cdot \x}\right) 
 -
\I \bk_0\cdot \x \,\E^{-\I \bk_i \cdot \x}\,\left( \Delta + |\bk_0|^2\right) \E^{\I \bk_j \cdot \x}\right)    \D V \nonu \\[2mm]
 &+ 2(\bk_0 \cdot \bk_0) \int_{\Pi} \E^{\I (\bk_j -\bk_i) \cdot \x} \, \D V=
 -2(\bk_i \cdot \bk_0)  \int_{\Pi}\E^{\I ( \bk_j - \bk_i) \cdot \x}\, \D V
 +2(\bk_0 \cdot \bk_0) \int_{\Pi} \E^{\I (\bk_j -\bk_i) \cdot \x} \, \D V \nonu \\[2mm]
 &= 2(\bk_0 - \bk_i) \cdot \bk_0\int_{\Pi} \E^{\I (\bk_j -\bk_i) \cdot \x} \, \D V =
 \left\{
 \begin{array}{cl}
 2(\mb_0 \cdot \bk_0) |\Pi|, & i = j =1, \\[2mm]
 0, & \text{otherwise}.
 \end{array}
 \right.
 \label{dv}
\end{align}
Thus, the latter formula with equality $ 2\mb_0 \cdot \bk_0 = |\mb_0|^2$ and \rf{Lij1}  implies \rf{L1}.
\end{proof}
\renewenvironment{proof}{{\noin \it Proof.}}{\qed}

\section{Local and global gaps in the Dirichlet problem}
\label{bgd}
\setcounter{equation}{0}

To find bands and gaps, it is necessary to construct a dispersion surface which we view as the graph of a periodic infinite-valued function of the Bloch vector $\bk$ in the entire space. The arguments at the beginning of \S\ref{np} imply that it is enough to study this surface only in a neighborhood of the cone $C_0: ~\om=c|\bk|$.

The following statement justifies the item (a) from the introduction.
\begin{theorem}
\label{5.1}
If $(\om_0, \bk_0)$ belongs to the cone $C_0$ but not to the intersection of the cones of the dispersion surface
of the unperturbed problem (see Figure \ref{cones}), then $g(\bk_0) = \emptyset$,
i. e. there is no local gap for wave vector $\bk_0$ and frequencies $\om$ near $\om_0$.
\end{theorem}
\begin{proof}
We need to show that local gap $g(\bk_0)$ cannot contain points $\om=\om_a$ such that $\om_a\to\om_0=c|\bk_0|$ as $a\to0$.
Let us consider the segment $\om=\om_a, ~\bk=(1+\delta)\bk_0$ in $\R^4$, where $|\delta|<\delta_0$ with some $a$-independent $\delta_0>0$. This segment intersects the cone $C_0$  when $\delta=\delta_a$ where $\delta_a\to0$ as $a\to0$. The same remains true with a different $\delta_a=\delta_a'$ if  $C_0$ is replaced by the dispersion surface of the problem with inclusions since the latter surface depends smoothly on $\bk$ and $a$ when points $\bk$ are not exceptional and $a$ is small, see Lemma \ref{l22}. Hence there is a Bloch wave with parameters $(\om_a, (1+\delta_a')\bk_0), ~ \delta_a' \to 0$
as $a \to 0$, and therefore $\om_a \notin g(\bk_0)$.
\end{proof}

We fix an interval $(\epsilon,\epsilon^{-1}$) of frequencies $\om$ with an arbitrary small $\epsilon>0$.

To prove item (b) we need the following Lemma concerning the
dispersion surface located above the interval $\bk=(1+\delta) \bk_0, ~|\delta| < \delta_0$ with some $a$-independent $\delta_0>0$.
For simplicity of the formulas below, we introduce in \rf{tau} new variables
\begin{align}
\ta=\dfrac{4\pi a q}{|\Pi|},\quad   {\tdelta}=\frac{\delta |\mb_0|^2}{2}.
\label{ad}
\end{align}
\begin{lemma}\label{L51}
 There is $\delta_0>0$ such that the dispersion surface in the neighborhood of the point $(\om_0, \bk_0) \in C_0,~\om_0 = c|\bk_0|,$ above the interval $\bk=(1+\delta) \bk_0$, $|\delta| < |\delta_0|$, is split  into the two branches
determined by
\begin{align}
\label{kpm}
 \om_{\pm}/c  &= |\bk_0| + \frac{1}{2|\bk_0|} \left( {\ta} + \nu {\tdelta} \pm \sqrt{{\ta}^2 + {\tdelta}^2} \right)
 + \O(\ta^2 + \tdelta^2),\quad  |\tdelta|\leq \delta_0, ~~\ta\ll1,
 \end{align}
 where
 \begin{align}
 \nu&= \frac{4|\bk_0|^2}{|\mb_0|^2} - 1 \geq 0.
 \label{nu}
\end{align}
\end{lemma}
\begin{proof}
 The dispersion relation near the point $\bk_0$ is determined by the equation $k=(1+\ep)|\bk|$, where $\ep=\ep(a,\bk)$  can be obtained from the EZE condition for the operator $\N:=\Np_{\bk,\varepsilon} - \Nm_{a,\varepsilon}$,
 see Lemma \ref{l21}.
 We studied this EZE problem in \S\ref{np} when $\bk=\bk_0$ by using the matrix representation \rf{conc6} of the operator $\N$, and this representation allowed us to reduce the calculation of the function $\ep=\ep(a,\bk)$ to the equation
 $\det \widetilde{\M} =0,$ where $\widetilde{\M}$ is a $2\times 2$ matrix that differs from the matrix $\M$  given by \rf{M1} only by the values of the remainder terms. The smoothness of the remainder terms of matrices $\M$ and
 $\widetilde{\M}$ was established using Lemma \ref{l21}.

 Theorem \ref{t41} provides  analogues of the matrix representation \rf{conc6} and formula \rf{M1} on the whole interval $\bk=(1+\delta) \bk_0$, $|\delta| < |\delta_0|$.
 The same arguments used to find $\ep$ for $\bk = \bk_0$ in \S\ref{np}, now allow us to find $\ep$ on the whole interval. These arguments lead to
\begin{align}
\det\left[2\ep |\bk|^2 |\Pi| \Im - 4\pi a q \J
+ \delta  |\Pi| \left[
 \begin{array}{cc}
  0 & 0 \\[2mm]
  0 & |\mb_0|^2
 \end{array}
\right]
+ \O(a^2 + \delta^2 + \ep^2 ) \right]=0, \quad \bk=(1+\delta) \bk_0,
\label{tau00}
 \end{align}
 where $|\delta| < |\delta_0|$ and the remainder term is infinitely smooth.

In new variables \rf{ad} equation \rf{tau00} is equivalent to
 \begin{align}
  \det & \left[
  \begin{array}{cc}
  \ept - \ta  + \O(\ta^2 + \tdelta^2 + \ept^2 ) & -\ta + \O(\ta^2 + \tdelta^2 + \ept^2 ) \nonu \\[2mm]
  -\ta + \O(\ta^2 + \tdelta^2 + \ept^2 ) & \ept - \ta + 2\tdelta + \O(\ta^2 + \tdelta^2 + \ept^2 )
  \end{array}
  \right] \\[2mm]
  &= \ept^2 +2(\tdelta - \ta)\ept - 2\tdelta \ta +\O(|\ta|^3 + |\tdelta|^3 + |\ept|^3 )=0,
  \label{tau0}
 \end{align}
where $\ept = 2\ep |\bk|^2$. If the term $|\ept|^3$ in the remainder is omitted, then the roots of the simplified equation \rf{tau0} are
\begin{align}
 \ept = \ta -\tdelta \pm \sqrt{\ta^2 + \tdelta^2} + \O(\ta^2 + \tdelta^2)
 \label{sim}
\end{align}
The same is true without the simplification of the equation but there is a difficulty to justify this fact because the equation degenerates at the point of interest  $\ta=\tdelta=0$, and we cannot apply the implicit function theorem to determine $\ept$. The same difficulty is present in equation \rf{tau00} if we view  $2\ep |\bk|^2 |\Pi|$ as the eigenvalue of the remaining terms of
the matrix \rf{tau00}. Then the zero eigenvalue of the matrix has multiplicity two when $\ta = \tdelta=0$.
Therefore we cannot guarantee a smooth dependence of the eigenvalue on the parameters. In our case the situation
is even worse since the remainder depends on the eigenvalue.

We will proceed as follows. We rewrite \rf{tau0} in the form
\begin{align}
 (\ept - \ta + \tdelta)^2 = \ta^2 +  \tdelta^2 +\O(|\ta|^3 + |\tdelta|^3 + |\ept|^3),
 \label{tau01}
\end{align}
divide this equation by $\ta^2 + \tdelta^2$, introduce polar coordinates $(\rho, \phi)$ in the plane $(\ta,\tdelta)$:
\begin{align}
 \rho = \sqrt{\ta^2 + \tdelta^2}, \quad \cos \phi = \ta/\rho, \quad \sin \phi = \tdelta/\rho,
\end{align}
and replace $\ept$ by $\tilde{\eta} = \ept / \sqrt{\ta^2 + \tdelta^2}$. In the new variables \rf{tau01} reads
\begin{align}
 (\tilde{\eta} - \cos \phi + \sin \phi)^2 = 1 +\O(|\ta| + |\tdelta| + \rho |\tilde{\eta}|).
 \label{tau02}
\end{align}
Now we can apply the implicit function theorem and solve \rf{tau02} for $\tilde{\eta}$. This leads to \rf{sim}. Since $\ept=2\ep |\bk|^2$ and $\ep |\bk|$ here can be replaced by $\frac{\om}{c}-|\bk|$, see \rf{k}, equation \rf{sim} can be rewritten as
\begin{align*}
\frac{\om}{c}-|\bk| = \frac{1}{2|\bk|}[\ta -\tdelta \pm \sqrt{\ta^2 + \tdelta^2} + \O(\ta^2 + \tdelta^2)].
\end{align*}
To complete the proof of the lemma, it remains to use the relation between $\bk$ and  $\bk_0$:
\[
\bk = (1+\delta) \bk_0 = (1 + 2\tdelta/|\mb_0|^2)\bk_0.
\]
\end{proof}

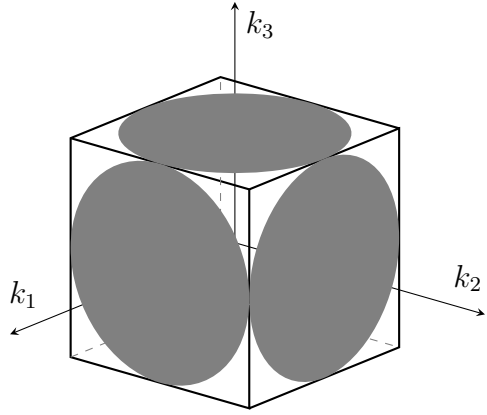
\begin{figure}[ht]
\centering
\begin{tikzpicture}[scale=1,>=stealth]
\begin{axis}[
    axis lines = center,
    xlabel = {$k_1$},
    ylabel = {$k_2$},
    zlabel = {$k_3$},
    xtick={-0.5,0.5},
    xticklabels={$-1/2$,$0.5$},
    yticklabels={$-1/2$,$0.5$},
    zticklabels={$-1/2$,$\frac{1}{2}$},
    ytick={-0.5,0.5},
    ztick={-0.5,0,0.5},
    xmin=0, xmax=1.5,
    ymin=0, ymax=1.4,
    zmin=-0.05, zmax=1.1,
    unit vector ratio = 1 1 1,
    view={130}{20}
]

  \draw[thick](0.5,0.5,-0.5)--(-0.5,0.5,-0.5)--(-0.5,0.5,0.5)--(0.5,0.5,0.5)--(0.5,0.5,-0.5)--(0.5,-0.5,-0.5)--(0.5,-0.5,0.5)--(-0.5,-0.5,0.5)--(-0.5,0.5,0.5);
  \draw[thick](0.5,0.5,0.5)--(0.5,-0.5,0.5);
  \draw[dashed, gray](0.5,-0.5,-0.5)--(-0.5,-0.5,-0.5)--(-0.5,0.5,-0.5);
  \draw[dashed, gray](-0.5,-0.5,-0.5)--(-0.5,-0.5,0.5);
  \node[right] at (0,0,0.5) {$0.5$};

\begin{scope}[canvas is xz plane at y=0.5]
    \fill[gray, opacity = 0.6] (0, 0) ellipse (0.5cm and 0.5cm);
\end{scope}

\begin{scope}[canvas is xy plane at z=0.5]
    \fill[gray, opacity = 0.6] (0, 0) ellipse (0.5cm and 0.5cm);
\end{scope}

\begin{scope}[canvas is yz plane at x=0.5]
    \fill[gray, opacity = 0.6] (0, 0) ellipse (0.5cm and 0.5cm);
\end{scope}
\end{axis}
\end{tikzpicture}
\label{gaps}
\caption{Shaded unit discs on the faces of the first Brillouin zone show the location of the Bloch vector
where waves cannot propagate.}
\end{figure}

\begin{theorem}
\label{5.2}
Let $\bk_0$ be an exceptional point of order two and $(\bk_0, \mb_0)$ satisfy \rf{Ew}.
\begin{enumerate}
\item
If $\dst \frac{|\bk_0|}{|\mb_0|} < \frac{\sqrt{2}}{2}$, then the local
gap $g(\bk_0)$ exists in a neighborhood of $\om_0=c|\bk_0|$ for small enough $a$ and  consists of frequencies $\om = \om_a$ for which
\begin{align}\label{int}
\dst |\bk_0|+\ta\nu_-+\O(a^2)< \om/c < |\bk_0|+\ta\nu_+ +\O(a^2), \quad \nu_\pm=1\pm \sqrt{1-\nu^2},
\end{align}
where $\ta$ is given by \rf{ad} and $\nu$ is given by \rf{nu}.
\item
If  $\dst \frac{|\bk_0|}{|\mb_0|} > \frac{\sqrt{2}}{2}$ and $a$ is small enough, then $g(\bk_0) \neq \emptyset$.
\end{enumerate}
\end{theorem}
\begin{remark}
The location of the Bloch vector on the surface of the first Brillouin zone for which waves
cannot propagate is shown in Figure \ref{gaps} by shaded unit disks.
\end{remark}

\begin{proof}
Local gap $g(\bk_0)$ consists of frequencies $\om=\om_a$ such that $\om_a\to c|\bk_0|$ as $a\to 0$ and the interval $l\in \mathbb R^4$ defined by equations $\om=\om_a, ~\bk=(1+\tdelta)\bk_0,~|\tdelta|<\tdelta_0,$ with some $a$-independent value of $\tdelta_0>0$ does not intersect the branches of the dispersion surface described in Lemma \rf{L51}.

We start the proof with the second part. Let $\dst \frac{|\bk_0|}{|\mb_0|} > \frac{\sqrt{2}}{2}$. Then for $\nu$ defined in Lemma \rf{L51} we have $\nu>1$. Using the inequality $\sqrt{\ta^2 + \tdelta^2}>|\tdelta|$, we obtain the following estimates for function $\om=\om_+$ defined in \rf{kpm} with small, $a$-independent $\tdelta_0>0$ and $a\to 0$:
\[
\left.\frac{\om_+}{c}\right|_{\tdelta=\tdelta_0}-|\bk_0|\geq(\nu+1)\tdelta_0/2>0, \quad \left.\frac{\om_+}{c}\right|_{\tdelta=-\tdelta_0}-|\bk_0|\leq(-\nu+1)\tdelta_0/2<0 .
\]
Hence, the range of the function $\om_+/c-|\bk_0|,~|\tdelta|<\tdelta_0,$ contains the segment $[(-\nu+1)\tdelta_0/2,(\nu+1)\tdelta_0/2]$ if $a$ is small enough. Thus, for any $\om_a$ close enough to $c|\bk_0|$, there exists $\tdelta_a$ such that $|\tdelta_a|<\tdelta_0$ and $\om_+=\om_a$. In other words, the ray $l$ intersect the graph of the dispersion surface $\om_+$ when $\tdelta=\tdelta_a,~\om_+=\om_a$, the Bloch wave propagates with the frequency $\om_a$ and $\bk=(1+\tdelta_a)|\bk_0|$,  and therefore $g(\bk_0) \neq \emptyset$. Note that the same result can be obtained using $\om_-$ instead of $\om_+$.

Let now $\dst \frac{|\bk_0|}{|\mb_0|} <\frac{\sqrt{2}}{2}$, and therefore $\nu<1$. Denote by $\om_\pm^0$ functions $\om_\pm$ defined in \rf{kpm} with the remainder terms omitted. The minimum of $\om_{+}^0/c-|\bk_0|$ is attained at $\tdelta = -\frac{\ta\nu}{\sqrt{1-\nu^2}}$ with
 the value $\om_{+}^0/c-|\bk_0| = \ta\nu_+$. Similarly, the maximum of $\om_{-}^0/c-|\bk_0|$ is attained at $\tdelta = \frac{\ta\nu}{\sqrt{1-\nu^2}}$ with
 the value $\om_{-}^0/c-|\bk_0| = \ta\nu_-$. Thus, the gap between the ranges of these two functions is the interval $(\ta\nu_-,~\ta\nu_+)$. If $|\tdelta|\leq C\ta$ with a large, but fixed $a$-independent constant $C$, then \rf{kpm} implies that the gap between ranges of functions $\om_{\pm}/c-|\bk_0|$ is the interval $I:=(\ta\nu_-+O(a^2),~\ta\nu_++O(a^2))$. When $C\ta\leq|\tdelta|\leq\tdelta_0$ and $\tdelta_0$ is small enough, the gap between functions $\om_{\pm}/c-|\bk_0|$ can be estimated through $\tdelta|$ and this gap includes interval $I$. Hence, $I$ is the exact interval such that the segment $l$ introduced at the beginning of the proof does not intersect the dispersion surface when $\om_a/c-|\bk_0|\in I$. This completes the proof of the theorem.
\end{proof}

Finally, we prove the theorem on the absence of global gaps under small perturbations.
\begin{theorem}
Global gaps do not exist in any fixed interval $\epsilon < \om < \epsilon^{-1}$ of the time frequency
$\om$ if the size $a$ of the inclusion is sufficiently small.
\end{theorem}
\begin{proof}
We will provide two proofs: the first one is geometrical and the second is shorter but more formal.
 For clarity, we explain the first proof in the two-dimensional case shown in Figure \ref{cones}.
 
 For $\om_0$ to be in a global gap, the plane $\om = \om_0$ must not intersect the dispersion
 surface of the problem with inclusions. Hence the surface in Figure \ref{cones} must break when $a>0$ and create a gap allowing 
 the plane $\om = \om_0$ to go through without intersection with the surface. It could happen only near the intersection of $C_0$ with other cones where the dispersion surface could split and move apart when $a > 0$.
 This splitting takes place, and the line of intersection of the cones where the splitting occurs is shown in Figure \ref{cones} by the dashed curve. Its projection on $\R^2_{\bk}$ is the line $AB$ consisting of exceptional points 
 of multiplicity two given by \rf{Ew1}. However, the plane $\om = \om_0$ intersects the cone $C_0$ along 
 the circle $\om_0 = c|\bk|$ and the splitting when $\bk$ is on the circle may occur only in a small neighborhood 
 of exceptional points near the circle. The plane intersects the dispersion surface in the neighborhood of other points on the circle and therefore $\om_0$ does not belong to a global gap.
 
An alternative proof is based on Theorem \ref{5.1}. The plane $\om = \om_0$ interests the cone $C_0$
 at the sphere $|\bk|^2 = \om_0^2/c^2$ which contains mostly non-exceptional values of $\bk$. For such non-exceptional point $\bk = \bk_0$, the set $g(\bk_0) = \emptyset$. According to Theorem \ref{5.1}, $\om_0$ cannot belong to the local gap $g(\bk_0)$ and therefore to a global as well.
\end{proof}

\section{Local and global gaps in the transmission problem}
\label{trans}
\setcounter{equation}{0}

The amplitude $u$ of the Bloch solution of the transmission problem satisfies the equation
\begin{align}
 \Delta u + k^2_{\pm} u= 0, \quad \x \in \R\setminus\de \Om,
 \label{}
\end{align}
with the boundary conditions on $\de \Om$
\begin{align}
 \llbracket  u(\x) \rrbracket =0, \quad \left\llbracket \frac{1}{\rho(\x)} \frac{\de u(\x)}{\de \n} \right\rrbracket =0,
\end{align}
where $k_{\pm} = \om/c_{\pm}$, $c_{\pm} = 1/\sqrt{\gamma_{\pm} \rho_{\pm}}$ is the speed of the wave propagation
in the medium and inclusion, respectively, $\gamma_{\pm}$ is the adiabatic bulk compressibility modulus, and $\rho(\x)$ is the mass density.
Here and on the subscript $\pm$ refers to the value of the quantity
outside/inside of the inclusion. The brackets $\llbracket  \cdot \rrbracket$ denote the jump of the enclosed quantity across $\de \Om$.

The analysis of gaps in the transmission problem is similar to that in the Dirichlet problem. In particular,
Lemma \ref{4.1} about matrix representation of $\Np_{\bk,\ep} - \Nm_{0,\ep}$  remains valid with the same representation
of matrix $\L$. Also, Theorem \ref{5.1} about the absence of local gaps near simple points $(\om_0, \bk_0)$ of the cone
$C_0$ still holds for the transmission problem.

While the exterior operators of the Dirichlet and transmission problems coincide, the interior operators have a different form. For the exceptional wave vector $\bk_0$ of order $n$ such that its integer-valued shifts $\bk_i = \bk_0 - \mb_i$,
$\mb_i \in \Z^3, ~0 \leq i \leq n-1$ are also satisfy \rf{u}, \rf{B1}, the asymptotics of the quadratic form
\begin{align}
 \left(\left(\Nm_{a,\ep} - \Nm_{0,\ep}\right) \psi_i, \psi_j \right), \quad \psi_j = \E^{-\I \bk_{j} \cdot \x}, \quad 0 \leq i,j \leq n-1,
\end{align}
is given by the matrix $\Mi$ obtained in \cite{GV:22}. For simplicity, we provide the form of this matrix in the case of spherical inclusions
 \begin{align}
M_{i,j} = |\Pi| |\bk_0^2| \left(\alpha + \beta \,\bkh_i \cdot \bkh_j \right) f + \O\left( a^4 + a^3 |\ep| \right), \quad 0 \leq i,j \leq n-1,
 \label{M}
\end{align}
where
\begin{align}
 \alpha = 1 - \frac{\gi}{\ge}, \quad
 \beta = 3\,\frac{\sigma - 1}{\sigma + 2}, \quad \sigma = \frac{\rho_{+}}{\rho_{-}}.
\end{align}
Here $\bkh = \bk/|\bk|$ and $f$ is the volume fraction of the inclusions. Matrix $\Mi$ plays the role of the matrix $4\pi a q \J$ in \rf{tau} for the Dirichlet problem. Therefore, if $\bk_0$ is an exceptional vector of multiplicity two, the $2\times 2$ matrix $\M$ in \rf{tau} has the form
\begin{align}
\M&=2\ep |\bk_0|^2 |\Pi| \Im
- f |\bk_0|^2  |\Pi| \left[
 \begin{array}{cc}
  \alpha + \beta & \alpha + \beta\,\bkh_0 \cdot \bkh_1\\[2mm]
  \alpha + \beta \,\bkh_0 \cdot \bkh_1& \alpha+ \beta
 \end{array}
\right]
+ \delta  |\Pi| \left[
 \begin{array}{cc}
  0 & 0 \\[2mm]
  0 & |\mb_0|^2
 \end{array}
\right] \nonu \\[2mm]
&+ \O(a^4 + a^3 |\ep|+ \delta^2 + \ep^2 ).
 \label{S}
\end{align}
Equating the leading term of the determinant of $\M$ to zero we obtain a quadratic equation for $\ep$
\begin{align}
 &4|\bk_0|^2 \ep^2 + 4\ep \left( \tdelta - f |\bk_0|^2 (\alpha + \beta) \right) 
 + f^2 |\bk_0|^2 \left( (\alpha + \beta)^2 - (\alpha + \beta \,\bkh_0 \cdot \bkh_1)^2 \right) -2\tdelta f (\alpha + \beta) =0.
\end{align}
The analog of Lemma \ref{L51} is
\begin{lemma}
 \label{L6.1}
 There is $\delta_0>0$ such that the dispersion surface in the neighborhood of the point $(\om_0, \bk_0) \in C_0,~\om_0 = c|\bk_0|,$ above the interval $\bk=(1+\delta) \bk_0$, $|\delta| < |\delta_0|$, is split  into the two branches
determined by
\begin{align}
\label{omtp}
\om_{\pm}/c &= |\bk_0| \left(1 + \oh (\alpha + \beta)f \right) \nonu \\[2mm]
&+ \frac{1}{2|\bk_0|}
\left(\nu \tdelta  \pm  \sqrt{(\alpha + \beta\,\bkh_0 \cdot \bkh_1)^2 |\bk_0|^4 f^2   + \tdelta^2} \right) + \O(a^4 +  \tdelta^2),\quad  |\tdelta|\leq \delta_0, ~~a\ll1,
 \end{align}
\end{lemma}
\noin
where $\tdelta$ and $\nu$ are given by \rf{ad} and \rf{nu}.
This expression is similar to \rf{kpm} and we can use the same analysis to find the location of local gaps.
The result is given by
\begin{theorem}
\label{T6.1}
Let $\bk_0$ be an exceptional point of order two and $(\bk_0,\mb_0)$ satisfy \rf{Ew}.
\begin{enumerate}
\item
If $\dst \frac{|\bk_0|}{|\mb_0|} < \frac{\sqrt{2}}{2}$, then the local
gap $g(\bk_0)$ exists in a neighborhood of $\om_0=c|\bk_0|$ for small enough $a$ and  consists of frequencies $\om = \om_a$ for which
\begin{align}
\label{}
|\tilde{\bk}_0| - \frac{\mu}{2|\bk_0|}\sqrt{1-\nu^2} + \O(a^4) < \om/c < |\tilde{\bk}_0| + \frac{\mu}{2|\bk_0|}\sqrt{1-\nu^2} + \O(a^4 ),
\end{align}
where $|\tilde{\bk}_0| = |\bk_0| \left(1 + \oh (\alpha + \beta)f \right)$, $\mu = |\alpha + \beta\,\bkh_0 \cdot \bkh_1| |\bk_0|^2 f$, and $\nu$ is given by \rf{nu}.
\item
If  $\dst \frac{|\bk_0|}{|\mb_0|} > \frac{\sqrt{2}}{2}$ and $a$ is small enough, then $g(\bk_0) \neq \emptyset$.
\end{enumerate}
\end{theorem}
Note that the existence of absolute gaps has been shown numerically for cubic arrays of air bubbles in water \cite{Kafesaki:00} but no absolute gaps were found for cubic arrays of rigid spheres in air \cite{Kushwaha:98a}.

\section{Conclusions}
\label{}
\setcounter{equation}{0}

We considered the propagation of acoustic waves in a medium containing a simple cubic array of small inclusions
of arbitrary shape. We show that global gaps do not exist in each fixed interval of the time frequency if the size of inclusions is small enough. The notion of local gaps defied by wave vectors $\bk$ is introduced and studied. We calculate the location of local gaps analytically for the Dirichlet and transmission problems. 
The width of the local bandgaps in the Dirichlet problem is proportional to the size of the inclusions while in 
the transmission problem problem, it is proportional to the volume fraction of the inclusions.
Although the bandgaps in the two problems were different, they exist for the same values of the Bloch vector $\bk$.
In particular,  local gaps exist on the boundary of the first Brillouin zone which is the unit cube in our scaling if $\bk$ belongs to the unit disks on the faces of the cube.

 \newcommand{\noop}[1]{}


\begin{thebibliography}{99}

\bibitem{Kushwaha:94}
Kushwaha MS, Halevi P, Mart\'{\i}nez G, Dobrzynski L, Djafari-Rouhani B. 1994
  Theory of acoustic band structure of periodic elastic composites. {\em
  Physical Review B} \textbf{49}, 2313--2322.

\bibitem{Johnson:2001}
Johnson S, Joannopoulos J. 2001  Block-iterative frequency-domain methods for
  Maxwell's equations in a planewave basis. {\em Optics Express} \textbf{8},
  173--190.

\bibitem{Wu:2004}
Wu T, Huang Z, Lin S. 2004  Surface and bulk acoustic waves in two-dimensional
  phononic crystal consisting of materials with general anisotropy. {\em
  Physical Review B} \textbf{69}.

\bibitem{Fan:1996}
Fan S, Villeneuve P, Joannopoulos J. 1996  Large omnidirectional band gaps in
  metallodielectric photonic crystals. {\em Physical Review B} \textbf{54},
  11245--11251.

\bibitem{Tanaka:2000}
Tanaka Y, Tomoyasu Y, Tamura S. 2000  Band structure of acoustic waves in
  phononic lattices: Two-dimensional composites with large acoustic mismatch.
  {\em Physical Review B} \textbf{62}, 7387--7392.

\bibitem{Taflove:2000}
Taflove A, Hagness SC. 2005 {\em Computational Electrodynamics: The
  Finite-Difference Time-Domain Method}.
Artech House Publishers, 3rd edition.

\bibitem{Kuchment:99}
Axmann W, Kuchment P. 1999  An efficient finite element method for computing
  spectra of photonic and acoustic band-gap materials - I. Scalar case. {\em
  Journal of Computational Physics} \textbf{150}, 468--481.

\bibitem{Dobson:1999}
Dobson D. 1999  An efficient method for band structure calculations in 2D
  photonic crystals. {\em Journal of Computational Physics} \textbf{149},
  363--376.

\bibitem{Giani:2012}
Giani S, Graham IG. 2012  Adaptive finite element methods for computing band
  gaps in photonic crystals. {\em Numerische Mathematik} \textbf{121}, 31--64.

\bibitem{JJWM:11}
Joannopoulos JD, Johnson SG, Winn JN, Meade RD. 2011 {\em Photonic Crystals:
  Molding the Flow of Light}.
Princeton, NJ: Princeton University Press.

\bibitem{Li:2013}
Li FL, Wang YS, Zhang C, Yu GL. 2013  Bandgap calculations of two-dimensional
  solid-fluid phononic crystals with the boundary element method. {\em Wave
  Motion} \textbf{50}, 525--541.

\bibitem{AM:76}
Ashcroft NW, Mermin ND. 1976 {\em Solid state physics}.
New York: Holt, Rinehart and Winston.

\bibitem{Kittel:04}
Kittel C. 2004 {\em Introduction to Solid State Physics}.
Wiley 8th edition.

\bibitem{Figotin:1996b}
Figotin A, Kuchment P. 1996  Band-gap structure of spectra of periodic
  dielectric and acoustic media. 2. Two-dimensional photonic crystals. {\em
  SIAM Journal on Applied Mathematics} \textbf{56}, 1561--1620.

\bibitem{Ammari:2009}
Ammari H, Kang H, Lee H. 2009  Asymptotic Analysis of High-Contrast Phononic
  Crystals and a Criterion for the Band-Gap Opening. {\em Archive for Rational
  Mechanics and Analysis} \textbf{193}, 679--714.

\bibitem{Lipton:2017}
Lipton R, Viator, Jr. R. 2017  Creating band gaps in periodic media. {\em
  Multiscale Modeling and Simulation} \textbf{15}, 1612--1650.

\bibitem{Lipton:2022a}
Lipton R, Viator, Jr. R, Bolanos SJ, Adili A. 2022  Bloch waves in high
  contrast electromagnetic crystals. {\em ESAIM-Mathematical Modelling and
  Numerical Analysis} \textbf{56}, 1483--1519.

\bibitem{GV:22}
Godin YA, Vainberg B. 2022  Clusters of Bloch waves in three-dimensional
  periodic media. {\em Proceedings of the Royal Society A: Mathematical,
  Physical and Engineering Sciences} \textbf{478}, 20220519.

\bibitem{GV:23}
Godin YA, Vainberg B. 2023  Propagation and dispersion of Bloch waves in
  periodic media with soft inclusions. eprint arXiv: 2207.08296.

\bibitem{Kafesaki:00}
Kafesaki M, Penciu RS, Economou EN. 2000  Air Bubbles in Water: A Strongly
  Multiple Scattering Medium for Acoustic Waves. {\em Phys. Rev. Lett.}
  \textbf{84}, 6050--6053.

\bibitem{Kushwaha:98a}
Kushwaha M, Djafari-Rouhani B, Dobrzynski L, Vasseur J. 1998  Sonic stop-bands
  for cubic arrays of rigid inclusions in air. {\em European Physical Journal
  B} \textbf{3}, 155--161.

\end{thebibliography}
\end{document}